\documentclass[12pt,a4paper]{amsart}

\makeatletter
\long\def\unmarkedfootnote#1{{\long\def\@makefntext##1{##1}\footnotetext{#1}}}
\makeatother

\usepackage{amsmath,amssymb}
\usepackage{enumerate,amssymb}
\usepackage{graphicx,enumerate}
\usepackage{color}
\theoremstyle{plain}
\newtheorem{thm}{Theorem}[section]

\newtoks\prt

\theoremstyle{definition}

\def\eqn#1$$#2$${\begin{equation}\label#1#2\end{equation}}

\numberwithin{equation}{section}

\def\epsilon{\varepsilon}

\def\er{\mathbb R}

\def\eps{\epsilon}

\def\Hom{{\it{Hom}}}
\def\Dif{{\it{Diff}}}

\def\loc{\operatorname{loc}}

\def\mir1{\mathcal L_1}

\def\phi{\varphi}

\def\rn{\mathbb R^n}

\newtoks\by
\newtoks\paper
\newtoks\book
\newtoks\jour
\newtoks\yr
\newtoks\pages
\newtoks\vol
\newtoks\publ
\def\ota{{\hbox\vol{???}}}
\def\cLear{\by=\ota\paper=\ota\book=\ota\jour=\ota\yr=\ota
\pages=\ota\vol=\ota\publ=\ota}
\def\endpaper{\the\by, {\the\paper},
\textit{\the\jour} \textbf{\the\vol} (\the\yr), \the\pages.\cLear}
\def\endbook{\the\by, \textit{\the\book}, \the\publ.\cLear}
\def\endprep{\the\by, \textit{\the\paper}, \the\jour.\cLear}
\def\endyearprep{\the\by, \textit{\the\paper}, \the\jour, (\the\yr).\cLear}
\def\name#1#2{#2 #1}
\def\nom{ \rm no. }

\title{Ball-Evans approximation problem: recent progress and open problems}

\author{Stanislav Hencl}
\address{Department of Mathematical Analysis, Charles University,
So\-ko\-lovsk\'a 83, 186~00 Prague 8, Czech Republic}
\email{\tt hencl@karlin.mff.cuni.cz}

\date{\today}

\newcounter{problem}

\begin{document}

\setcounter{problem}{0}
\def\open#1{\vskip 5pt\addtocounter{problem}{1} {\bf Open problem \theproblem.} {\underline{#1}}:}

\thanks{The author was supported by the grant GA\v{C}R P201/24-10505S}

\dedicatory{Dedicated to Vladimir Gol'dshtein.}

\begin{abstract}
    {In this paper we give a short overview about the Ball-Evans approximation problem, i.e. about the approximation of Sobolev homeomorphism by a sequence of diffeomorphisms (or piecewise affine homeomorphisms) and we recall the motivation for this problem. We show some recent planar results and counterexamples in higher dimension and we give a number of open problems connected to this problem and related fields.}
\end{abstract}

\maketitle

\section{Introduction and motivation}

In this paper we would like to recall some recent progress in the solution of the Ball-Evans approximation problem and to state some open problems in the field. 
In the whole paper $\Omega\subset\rn$, $n\geq 2$, denotes an open set and $f$ is a mapping from $\Omega$ to $\rn$. The Ball-Evans approximation problems asks for the approximation of homeomorphism $f$ in the Sobolev space $W^{1,p}(\Omega,\rn)$ by a sequence of diffeomorphisms (or piecewise affine homeomorphisms) in the norm of Sobolev space. 

Let us start with more basic and general problem of approximation of homeomorphisms $f:\Omega\to f(\Omega)$ with piecewise affine homeomorphisms in the $L^{\infty}$ norm which has a long history. In the simplest non-trivial planar ($n=2$) setting the problem was solved by Rad\'{o}~\cite{Rado}. Due to its fundamental importance in geometric topology, this problem for $n>2$ was deeply investigated in the '50s and '60s. In particular, it was solved by Moise~\cite{Moise1} and Bing~\cite{Bing} in the case $n=3$ (see also the survey book~\cite{Moise2}), while for contractible spaces of dimension $n\geq5$ the result follows from theorems of Connell~\cite{Conn}, Bing~\cite{Bing2}, Kirby~\cite{Kirby} and Kirby, Siebenmann and Wall~\cite{Kirbetal} (for a proof see, e.g., Rushing~\cite{Rush}). Finally, twenty years later Donaldson and Sullivan proved that the result is false in dimension 4: more precisely, there is a homeomorphism from the unit ball of $\er^4$ to $\er^4$ which cannot be approximated by bilipschitz homeomorphisms, see~\cite[Corollary, page~183]{DonSull}.\par

Once completely solved in the uniform sense, the approximation problem became of interest again for variational models in Nonlinear Elasticity. In the setting of nonlinear elasticity (see for instance the pioneering work by Ball~\cite{Ball1} or monograph of Ciarlet \cite{Ci}) we study existence and regularity properties of minimizers of energy functionals of the form 
\begin{equation}\label{energyI}
I(f)=\int_{\Omega} W(Df(x))\,dx\,,
\end{equation}
where $f:\rn\supseteq\Omega\to\Delta\subseteq\rn$ ($n=2,3$) models the deformation of a homogeneous elastic material with respect to a reference configuration $\Omega$ and  $W:\er^{n\times n}\to\er$ is the stored-energy functional. There are some natural assumptions about the model so that the deformations $f$ are physically relevant.  
As pointed out by Ball in~\cite{Ball,Ball2} the deformation should be continuous as the material cannot break during our deformation and the 
"non-interpenetration of the matter" requires that our deformation is $1$ to $1$ and thus we study deformations that are in fact a homeomorphisms. Moreover, we usually require that 
\eqn{E_WA}
$$
\begin{aligned}
W(A)&\to+\infty\text{ as }\|A\|\to \infty,\\ 
W(A)&\to+\infty\text{ as }\det A\to 0\text{ and }\\ 
W(A)&=+\infty\quad\text{if $\det A\leq0$}.\\
\end{aligned}
$$
The first condition in~\eqref{E_WA} tells us that it costs too much to stretch our material too far, the second one tells us that it costs us too much to squeeze too much and the last condition tells us that the orientation is preserved during our deformation. 
In fact for physically relevant energy functionals we have $W(Df(x))\geq |Df(x)|^p$ for some $p\geq 1$ and thus our deformation $f$ of finite energy naturally belongs to some first order Sobolev space $W^{1,p}(\Omega,\rn)$. 


As pointed out by Ball in~\cite{Ball, Ball2} (who ascribes the question to Evans~\cite{Evans}), an important issue toward the understanding of the regularity of the minimizers in this setting (i.e., $W$ quasiconvex and satisfying~\eqref{E_WA}) would be to show the existence of minimizing sequences given by piecewise affine homeomorphisms or by diffeomorphisms. 
The first step is to prove that any homeomorphism $f\in W^{1,p}(\Omega,\rn)$, $p\in[1,+\infty)$, can be approximated by diffeomorphisms (or by piecewise affine homeomorphisms) in the $W^{1,p}$ norm and this problem is nowadays called the Ball-Evans approximation problem. The main motivation is that the usual approach for proving regularity of minimizers is to test the weak equation or the variation formulation by the solution itself. 
Unfortunately,  the integrals in question do not need to be finite and hence we would like to test the equation with a smooth test mapping which is close to the original $f$. Of course this approximation has to be in the correct class of competitors, i.e. it has to be a homeomorphism. 
There are of course other applications of such a result as the usual proof in the theory of Sobolev mappings is first shown for smooth mappings by direct computation and by approximation it then follows for general mapping. Therefore the positive solution to Ball-Evans approximation problem  
would significantly simplify many other known proofs, and it would easily lead to stronger new results. 


\begin{figure}
\centering
\vskip -20pt
\includegraphics[scale=0.3]{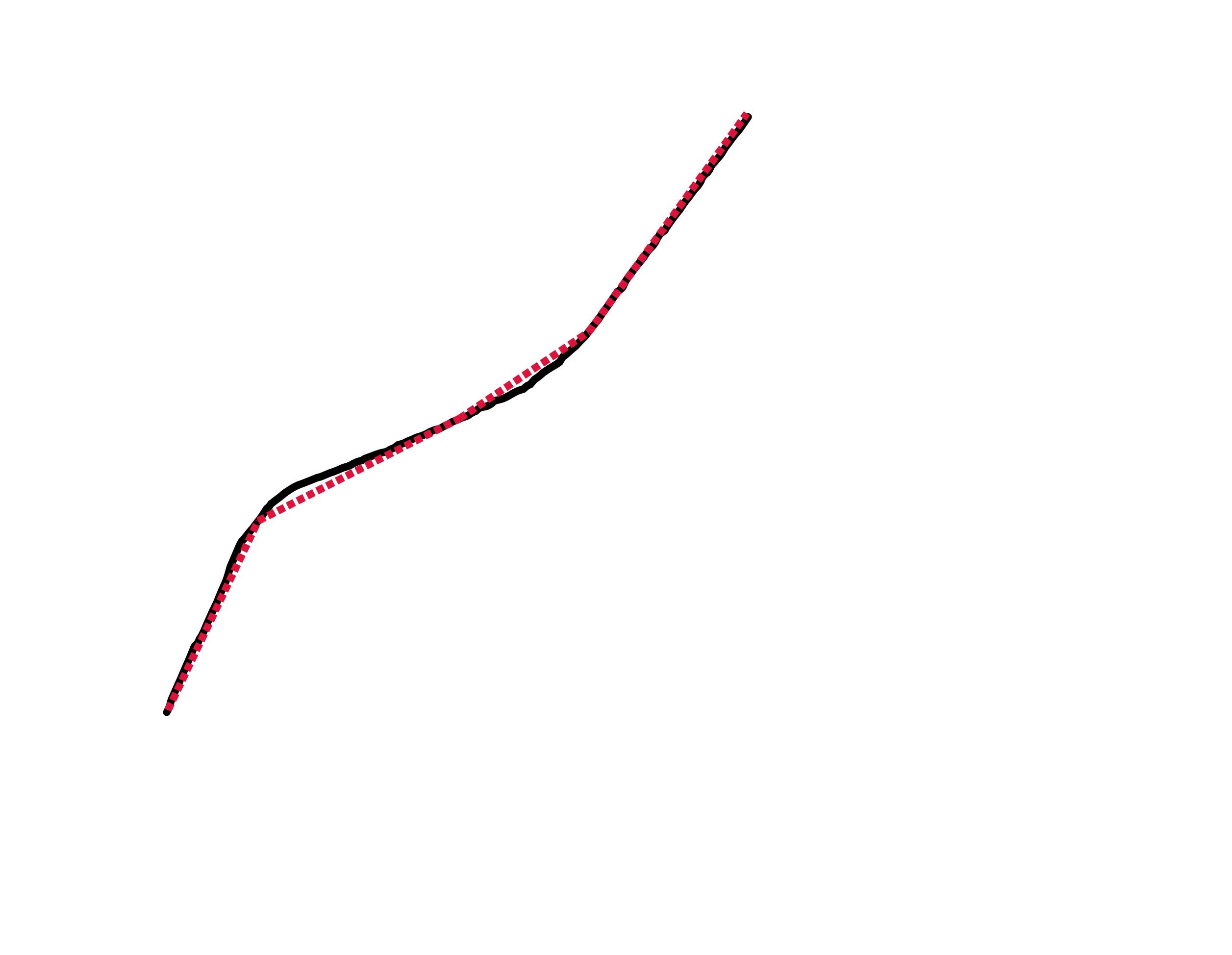}
\vskip -20pt
\caption{Approximation by piecewise linear homeomorphism is simple for $n=1$. We approximate black homeomorphism=increasing function by piecewise linear red one.}\label{fign1}
\end{figure}

\begin{figure}
\centering
\vskip -40pt
\includegraphics[scale=0.5]{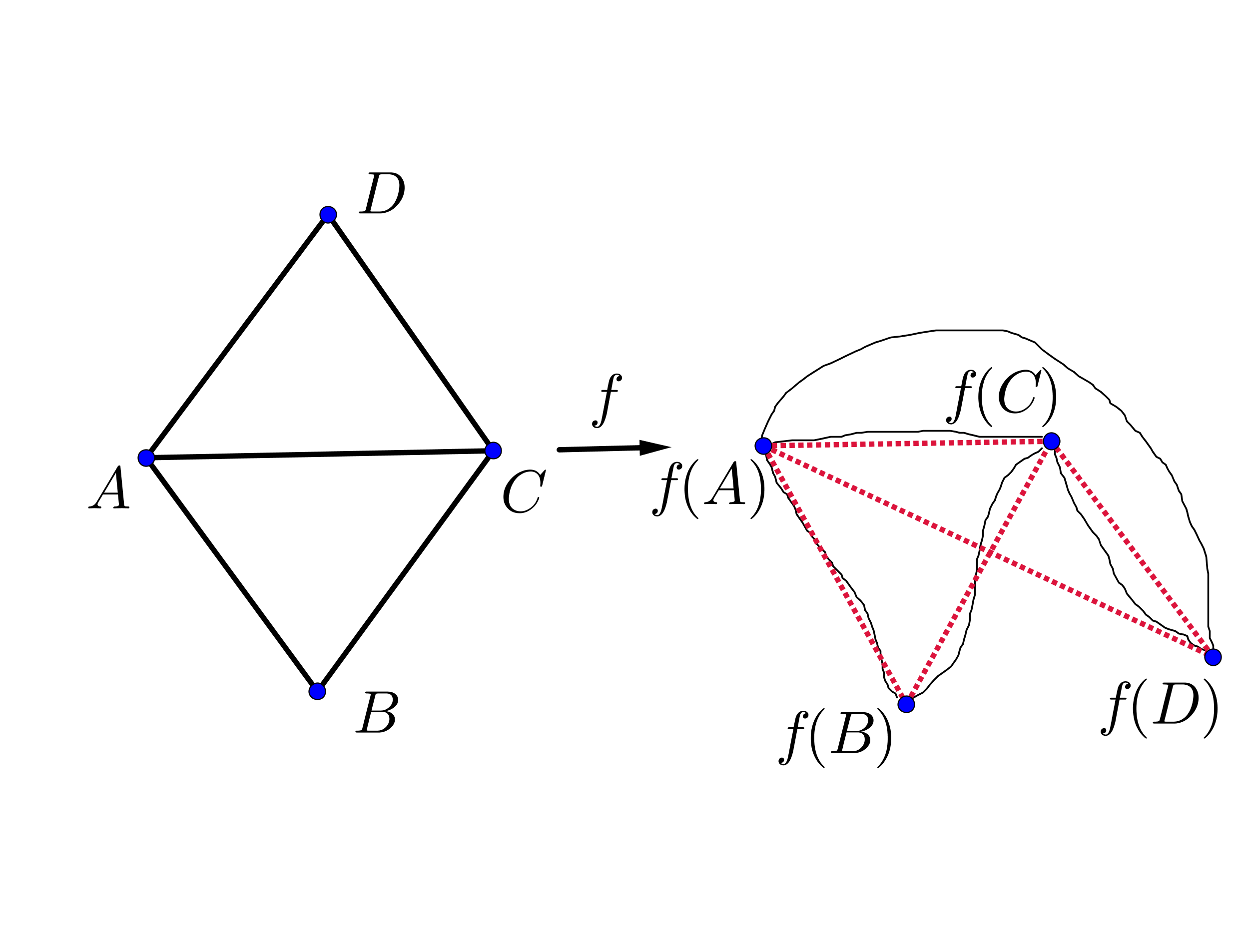}
\vskip -40pt
\caption{Approximation by piecewise linear homeomorphism is difficult already for $n=2$. The black curves show the image of corresponding segments by the real homeomorphism $f$. The natural piecewise linear approximation (denoted by red dotted segments) maps triangle $ABC$ onto triangle $f(A)f(B)f(C)$ and triangle $ACD$ onto triangle $f(A)f(C)f(D)$. We see that triangles $f(A)f(B)f(C)$ and $f(A)f(C)f(D)$ intersect and hence our approximation is not $1$ to $1$ there.}\label{fign2}
\end{figure}

It is easy to see that the approximation of homeomorphism by piecewise linear homeomorphisms (or diffeomorphisms) is simple in dimension $n=1$ (see Fig. \ref{fign1}) as each homeomorphism is if fact increasing (or decreasing) function. We can just choose any fine partition of the interval and connect the values by piecewise linear map which is again increasing and it is not difficult to show that the $\int_a^b |Df|^p$ is close enough for fine enough partition. 
However, the finding of diffeomorphisms near a given homeomorphism is not an easy task already in dimension $n=2$, as the usual approximation techniques like mollification or Lipschitz extension using maximal operator destroy, in general, the injectivity (see Fig. \ref{fign2}). Of course we need to approximate our homeomorphism $f$ not just with smooth maps, but with smooth homeomorphisms (otherwise the approximating sequence would be not even admissible for the original problem).


\section{Recent results}

Let us now describe the results which are known in the literature in this direction. The first ones were obtained in 2009 by Mora-Corral~\cite{M} (for planar bi-Lipschitz mappings that are smooth outside a finite set) and by Bellido and Mora-Corral~\cite{BMC}, in which they prove that if $u, u^{-1}\in C^{0,\alpha}$ for some $\alpha\in (0,1]$, then one can find piecewise affine approximations $v$ of $u$ in $C^{0,\beta}$, where $\beta\in(0,\alpha)$ depends only on $\alpha$.\par

The fundamental breakthrough result which strengthened the research interest in this area was obtained by Iwaniec, Kovalev and Onninen~\cite{IKO} (see also \cite{IKO2}) where they approximated planar $W^{1,p}$ homeomorphism for $1<p<\infty$. 
The remaining case $p=1$ was shown later by Hencl and Pratelli \cite{HP} by entirely different method and the method of \cite{HP} was later used e.g. by Campbell \cite{Ca} to obtain the approximation also in other function spaces. The main result of \cite{IKO2} and \cite{HP} can be formulated as follows. \par 

\begin{thm}\label{main}
Let $\Omega\subseteq\er^2$ be an open set, $1\leq p<\infty$ and $f\in W^{1,p}(\Omega,\er^2)$ be a homeomorphism. For every $\eps>0$ there is a smooth diffeomorphism (as well as a countably --but locally finitely-- piecewise affine homeomorphism) $f_\eps\in W^{1,p}(\Omega,\er^2)$ such that 
$$
\|f_\eps-f\|_{W^{1,p}}+\|f_\eps-f\|_{L^\infty}<\eps. 
$$
Moreover, $f_\eps(\Omega)=f(\Omega)$ and in particular, if $f$ is continuous up to the boundary of $\Omega$, then the same holds for $f_\eps$, and $f_\eps=f$ on $\partial\Omega$.
\end{thm}

 The main open problem in the area is to obtain any positive result in higher dimensions and especially in the physical relevant dimension $n=3$.

\open{Ball-Evans problem in  dimension $n\geq 3$} Prove or disprove that $W^{1,p}$ Sobolev homeomorphism in dimension $n\geq 3$ can be approximated by a sequence of diffeomorphism (or piecewise affine homeomorphisms). The positive result would be interesting for any value of $p$ and we could expect that the proof would be easier for $p>n-1$ (or even $p>n$) 
as we know many useful properties of Sobolev homeomorphisms for $p>n-1$ like differentiability a.e. (see e.g. \cite[Corollary 2.25]{HK}) and so on. It is also not clear at the moment if this assumption is needed and in principle there might be some counterexamples for small values of $p$. 

\vskip 10pt

In fact, we have two possible problems, one to approximate by diffeomorphisms and the second one to approximate by piecewise affine homeomorphism. Actually, both results would be interesting in different contexts. Approximation by smooth mappings would be useful for the possible proof of regularity of minimizers and the approximation by piecewise affine appears naturally in problems connected with finite elements methods in Numerical mathematics as it is a priori not clear if there even exists a finite element approximation which is $1$ to $1$.  
Luckily, the two things are equivalent at least in the planar case. 
Given an approximation by diffeomorphism it is possible to construct an piecewise affine approximation by homeomorphisms in any dimension as shown by Iwaniec and Onninen in \cite{IO3}. 
The converse is not immediate but, at least in the plane, it was shown by Mora-Corral and Pratelli in~\cite{MP}. 
It seems very reasonable that similar should be true at least for $n=3$ but it was not shown yet. 


\open{Piecewise affine implies diffeomorphic approximation for $n\geq 3$} Let $\Omega\subset\rn$ be open and $n\geq 3$. 
Show that given piecewise affine homeomorphism $f:\Omega\to\rn$ we can approximate it by diffeomorphisms in the Sobolev norm. 
It is quite likely that the possible solution to Open problem 1 will be done first for piecewise affine approximation (as in \cite{HP} for the planar case) and then we would need this result to get it also for diffeomorphisms. 

\vskip 10pt

The solution of Ball-Evans approximation in higher dimension might be much more difficult as the result is not true for small values of $p$ for $n\geq 4$. This was first shown by Hencl and Vejnar \cite{HV} for $p=1$ and later generalized for $1 \le p < [n/2]$ by Campbell, Hencl and Tengvall \cite{CHT}. Here, $[p]$ denotes the integer part of $p\in\er$, i.e. the smallest integer smaller or equal to $p$. 

\begin{thm}\label{CrazyMap}
Let $n \ge 4$ and $1 \le p < [n/2]$. Then there is a homeomorphism $f \in W^{1,p}((-1,1)^{n}, \rn)$ such that $J_{f}>0$ on a set of positive measure and $J_{f}<0$ on a set of positive measure. It follows that there are no diffeomorphisms (or piecewise affine homeomorphisms) $f_{k} : (-1,1)^{n} \to \rn$ such that $f_{k} \to f$ in $W_{\loc}^{1,p}((-1,1)^{n}, \rn)$. 
\end{thm}

It is not difficult to show that once we have a homeomorphism where the Jacobian is changing sign then we cannot approximate it. 
Indeed, assume on the contrary that $f$ from the statement can be approximated by diffeomorphisms (or piecewise affine homeomorphisms) $\{ f_{k} \}_{k=1}^{\infty}$, then the pointwise limit of a subsequence (which we denote the same) satisfies 
\begin{align*}
Df_{k}(x) \to Df(x) \quad \text{and} \quad J_{f_{k}}(x) \to J_{f}(x) 
\end{align*}
for almost every $x \in (-1,1)^{n}$. As $f_{k}$ are locally Lipschitz we know that $J_{f_{k}} \ge 0$ a.e. in $(-1,1)^{n}$ or $J_{f_{k}} \le 0$ a.e. in $(-1,1)^{n}$, see e.g. \cite{HM} and \cite[Theorem 5.22]{HK}. The pointwise limit of nonnegative (or non-positive) functions $J_{f_k}$ cannot change sign which gives us contradiction.

Let us note that the Jacobian of a Sobolev homeomorphism cannot change sign for $p>[n/2]$ as shown by Hencl and Mal\'y in \cite{HM} and thus the above result is sharp at least for this problem which was formulated by P. Hajlasz in 2001. However, it is not clear if the condition $p<[n/2]$ is sharp for the counterexamples to Ball-Evans approximation problem and there might be different counterexamples based on something entirely different.  

\open{Possible counterexample for $p\geq [n/2]$} Is it possible to construct a counterexample to Ball-Evans approximation problem for some $p\geq [n/2]$ in higher dimensions? 

\vskip 10pt

As we have already mentioned one of the main motivations for the study of Ball-Evans approximation problem is the possible application to the proof of regularity of minimizers of energy functional in Nonlinear Elasticity. However, if we want to really apply it there we need to be able to approximate all term in the energy functional. As the model examples we could mention for example the following two. 

\open{Ball-Evans problem for more general functionals} Given an open set $\Omega\subset\er^2$ and a homeomorphism $f:\Omega\to\er^2$ with $J_f>0$ a.e., $p\geq 1$ and $a>0$ let us define 
$$
E_1(f):=\int_{\Omega} \Bigl(|Df(x)|^p+\frac{1}{J_f^a(x)}\Bigr)\; dx\text{ and } E_2(f):=\int_{\Omega} \Bigl(|Df(x)|^2+\frac{|Df(x)|^2}{J_f(x)}\Bigr)\; dx.
$$
Can we find an approximating sequence of diffeomorphisms which approximate $f$ well both in the Sobolev norm and also in the $E_1$ or $E_2$ functional? 

\vskip 10pt

The second problem could be in fact equivalently rewritten by a simple change of variables as 
$$
E_2(f):=\int_{\Omega} \Bigl(|Df(x)|^2+\frac{|Df(x)|^2}{J_f(x)}\Bigr)\; dx=\int_{\Omega} |Df(x)|^2\; dx+\int_{f(\Omega)} |Df^{-1}(y)|^2\; dy
$$
and asks us for the approximation of bi-Sobolev homeomorphisms, i.e. to approximate both $f$ and $f^{-1}$ in the corresponding norms. This problem of bi-Sobolev approximation was first formulated in \cite{IKO} for any $W^{1,p}$, $p\geq 1$. So far it was solved only for $n=2$ and $p=1$ by Pratelli \cite{BiP} and for $n=2$ in general only under the additional assumption that $f$ is bi-Lipschitz by Daneri and Pratelli in~\cite{DP} and~\cite{DP2}. 

\vskip 10pt

In models on Nonlinear Elasticity we can work with compressible material (like rubber) or some incompressible one (like steel). 
Therefore in many modeling situations we actually work with incompressible deformations assuming that $J_f\equiv 1$ a.e. 
Naturally, we can ask if incompressible homeomorphisms in the Sobolev space could be approximated by incompressible diffeomorphisms and this could be important e.g. in problems related to \cite{AABT}. In fact the following open problem was asked to me by Assis Azavedo. 

\open{Ball-Evans mapping for incompressible homeomorphisms} Let $f\in W^{1,p}(\Omega,\rn)$ be a homeomorphism with $J_f\equiv 1$ a.e. Can we find a sequence of incompressible diffeomorphisms (or piecewise affine incompressible homeomorphisms) which approximate $f$ in the Sobolev norm? Note, that this is open already in the planar case $n=2$. 

\vskip 10pt

Some models of Nonlinear Elasticity include in fact also the second gradient and we work with the functional 
	\eqn{functional2}
	$$
	E(f)=\int_{\Omega}\bigl(W(Df(x))+\delta_0|D^2f(x)|^q \bigr)\; dx,
	$$
	where $q\in[1,\infty)$ and $\delta_0>0$. These models  
	were introduced by Toupin \cite{T}, \cite{T2}
	and later considered by many other authors, see e.g. Ball, Curie, Olver \cite{BCO}, Ball, Mora-Corral \cite{BMC}, M\"uller \cite[Section 6]{Mbook}, Ciarlet \cite[page 93]{Ci} and references given there. 
	The contribution of the higher gradient is usually connected with interfacial energies and is used to model various phenomena like elastoplasticity or damage.

\open{Ball-Evans approximation problem for $W^{2,q}$ homeomorphisms} Let $f\in W^{2,q}(\Omega,\rn)$, $g\geq 1$, be a homeomorphism. Can we find a sequence of diffeomorphisms (or piecewise affine homeomorphisms) which approximate $f$ in the $W^{2,q}$ Sobolev norm? Again this problem is open already in the planar case $n=2$ and it seems to be much more difficult than the corresponding $W^{1,p}$ problem. For some initial partial results in the $W^{2,1}$ case see \cite{CH}.

\vskip 10pt

\section{Weak and strong closure of homeomorphisms and diffeomorphisms}

Techniques developed in \cite{IKO} and \cite{HP} for the solution of the planar Ball-Evans approximation problem were later used for the solution of other problems. Once we apply standard techniques of Calculus of Variations to the study of minimizers of energy functionals we naturally obtain a sequence $f_k$ (in the class of energy competitors) whose energy $E(f_k)$ is converging to the infimum of the energy. By some natural coercivity we obtain that $f_k$ is actually converging to $f$ weakly in the corresponding Sobolev space and we would like to infer some properties of $f$ based on the properties of $f_k$. 

One can ask that the energy competitors are injective and regular and possible discontinuities happen only for the limit mapping $f$. 
In this way the possible deformation $f$ is reasonably close to injective mapping and corresponds to reasonable real-world deformation. 
Note that some phenomena like cavitation can really happen in some physical experiments like deformation of rubber cylinders (see e.g. references in \cite{Bc}) and this sort of discontinuity is easy to approximate by homeomorphisms. Thus we can ask that our class of competitors consists of homeomorphisms or even diffeomorphisms and then the minimizers of the energy belong to the class of weak limits of homeomorphisms (or weak limits of diffeomorphisms) from bounded open set $\Omega\subset\rn$ to bounded open set $\Omega'\subset\rn$ 
$$
\begin{aligned}
\overline{\Hom(\Omega,\Omega')}^{w}:=\Bigl\{&f\in W^{1,p}(\Omega,\Omega'):\ \text{ there are homeomorphisms }f_k\\
&\text{ from }\Omega\text{ onto }\Omega'\text{ with }f_k\to f\text{ weakly in }W^{1,p}\Bigr\},\\
\overline{\Dif(\Omega,\Omega')}^{w}:=\Bigl\{&f\in W^{1,p}(\Omega,\Omega'):\ \text{ there are diffeomorphisms }f_k\\
&\text{ from }\Omega\text{ onto }\Omega'\text{ with }f_k\to f\text{ weakly in }W^{1,p}\Bigr\}.\\
\end{aligned}
$$
This class of mapping in the plane $n=2$ was studied and characterizes by Iwaniec and Onninen in \cite{IO} and \cite{IO2} and by De Philippis and Pratelli in \cite{DPP} and they were able to prove the remarkable result that the class of weak limits or strong limits are actually the same in the planar case
$$
\overline{\Hom(\Omega,\Omega')}^{w}=\overline{\Hom(\Omega,\Omega')}^{\|\cdot\|_{W^{1,p}}}=\overline{\Dif(\Omega,\Omega')}^{w}=\overline{\Dif(\Omega,\Omega')}^{\|\cdot\|_{W^{1,p}}},
$$
where $\overline{\Hom(\Omega,\Omega')}^{\|\cdot\|_{W^{1,p}}}$ and $\overline{\Dif(\Omega,\Omega')}^{\|\cdot\|_{W^{1,p}}}$ denote the class of strong limits of homeomorphisms and diffeomorphisms in the $W^{1,p}$ norm. 

\open{Weak and strong closure of homeomorphisms for $n\geq 3$} Let $n\geq 3$, is it true that 
$$
\overline{\Hom(\Omega,\Omega')}^{w}=\overline{\Hom(\Omega,\Omega')}^{\|\cdot\|_{W^{1,p}}}?
$$
Somewhat surprisingly the result is false for $n=3$ and $p=2$ by the result of \cite{DHM} and it seems reasonable that this could be generalized for any $n\geq 3$ and $p=n-1$. However, this idea for counterexample could work only for $p=n-1$ and it is not clear at all what happens for other values of $p$. 

\vskip 10pt

\open{Weak and strong closure of diffeomorphisms for $n\geq 3$} 
Let $n\geq 3$, is it true that 
$$
\overline{\Dif(\Omega,\Omega')}^{w}=\overline{\Dif(\Omega,\Omega')}^{\|\cdot\|_{W^{1,p}}}?
$$

If we are able to solve Ball-Evans approximation problem for $n\geq 3$ then there is a chance that the method of the proof could be used for the positive solution of the two problems above (similarly to \cite{IO}, \cite{IO2} and \cite{DPP}). On the other hand, for small values of $p$ there might be other counterexamples based on a different idea than \cite{DHM}. For example the counterexample to Ball-Evans approximation problem for $n\geq 4$ and $1\leq p<[n/2]$ (see \cite{HV} and \cite{CHT}) cannot be obtained as a strong limit of diffeomorphisms (as Jacobians of diffeomorphisms do not change sign) but maybe it could be obtained as a weak limit of diffeomorphisms.

\vskip 10pt

\section{Can Jacobian change sign or even vanish}

In 2001 P. Hajlasz asked a problem if the Jacobian of a Sobolev homeomorphism can change sign (see e.g. the introduction in \cite{GH}). The positive solution for $p>[n/2]$ (or $p\geq 1$ for $n=3$) was given in \cite{HM} and  
as we have seen in Theorem \ref{CrazyMap} the counterexamples to Hajlasz problem for $1\leq p<[n/2]$ give also counterexamples to Ball-Evans approximation problem. In this section we would like to recall some problems remaining in this area. 

\vskip 10pt

\open{The limiting case $p=[\frac{n}{2}]$} Let $n\geq 4$ and $p=[\frac{n}{2}]$. Let $\Omega\subset\rn$ be connected open set and let $f\in W^{1,[\frac{n}{2}]}(\Omega,\rn)$ be a homeomorphism. Is it possible that both sets $\{J_f>0\}$ and $\{J_f<0\}$ have positive measure? For some partial result in this direction see Goldstein and Hajlasz \cite{GH1}. 

\vskip 10pt

In the original counterexamples in \cite{HV} and \cite{CHT} we have in fact a homeomorphisms which equals to identity far away from the origin (and thus the Jacobian is positive there) and its Jacobian is negative on some Cantor type set of positive measure and the set $\{J_f<0\}$ is thus totally disconnected. 
By some very careful and delicate iteration of similar construction we can in fact achieve (see \cite{CDH}) that $f$ equals to identity far away from the origin and $J_f<0$ a.e. in $[0,1]^4$ for any $p<\frac{3}{2}$. It is not clear if this can be achieved for all values $1\leq p<[n/2]$. 

\open{Can we have $J_f<0$ a.e. for all $p<[n/2]$} Let $n\geq 4$ and $1\leq p<[n/2]$. Is it possible to construct a homeomorphism $f\in W^{1,p}(\overline{B(0,R)},\overline{B(0,R)})$ such that $f(x)=x$ for $x\in \partial B(0,R)$ and $J_f<0$ a.e. in $[0,1]^n$? 

\vskip 10pt

The usual way how to prove that the mapping $f:\Omega\to\rn$ in question is in fact a homeomorphism is to first show that it is continuous, open (maps open sets to open sets) and discrete (preimage of single point does not have an accumulation point in $\Omega$). Then we prescribe a homeomorphic boundary condition on $\partial\Omega$ and it follows that $f$ is in fact a homeomorphism (see e.g. \cite[Chapter 3]{HK} for details). The question about sign of the Jacobian for continuous, open and discrete mappings seems to be more difficult. 

\open{Sign of the Jacobian for continuous, open and discrete mappings} Let $n\geq 3$ and $[n/2]<p\leq n-1$ (or $1\leq p\leq 2$ for $n=3$). Let $\Omega\subset\rn$ be connected open set and let $f\in W^{1,p}(\Omega,\rn)$ be a continuous, open and discrete mapping. Is it possible that both sets $\{J_f>0\}$ and $\{J_f<0\}$ have positive measure? 

Note that for $p>n-1$ the result should be true as these maps should be differentiable a.e. and then we could use the usual approach with the degree. Similarly to homeomorphisms each such a map should be either sense-preserving or sense-reversing, but we cannot use the theory of linking number (which was crucial in \cite{HM}) since the linking number might not be preserved under our mapping $f$. 

\vskip 10pt

At the end of Section 2 we have briefly discussed the models with second derivatives and their motivation. The problems discussed here could be easily asked also for second order Sobolev spaces.  

\open{Sign of the Jacobian for $W^{2,p}$ mappings} Let $n\geq 6$ and $1\leq q\leq \frac{[n/2]}{2}$. Let $\Omega\subset\rn$ be connected open set and let $f\in W^{2,q}(\Omega,\rn)$ be a homeomorphism. Is it possible that both sets $\{J_f>0\}$ and $\{J_f<0\}$ have positive measure? 

Of course, by the Sobolev embedding theorem we know that for $1\leq q<n$ we have $W^{2,q}\subset W^{1,q^*}$ for $q^*=\frac{nq}{n-q}$. If $q^*>[n/2]$ then the result must be true by the corresponding result for first order Sobolev space. 

In fact, there is a high chance that we can say even more. The proof in \cite{HM} is essentially based on the theory of linking number, i.e. topological concept which tells us how many times one curve (or manifold) goes around another one. In the proof we fix two cleverly chosen spheres of dimension $[n/2]$ (or one of dimension $[n/2]$ and one of dimension $[n/2]-1$ for even $n$) inside $\rn$. We consider the restriction of our $f$ on these spheres and we compute the linking number of the image of these spheres. We can assume that $f\in W^{1,p}$ on those spheres and the main argument is based on the fact that for $p>[n/2]$ our Sobolev space $W^{1,p}$ is embedded into the space of continuous functions on those $[n/2]$-dimensional spaces. Once we work with $W^{2,q}$ functions we can assume that our $f$ belong to $W^{2,q}$ on those spheres and for $q>\frac{[n/2]}{2}$ we know that $W^{2,q}$ is embedded into continuous functions (for $n=4,5$ we have $[n/2]=2$ and $W^{2,1}$ is also embedded into continuous functions on $2$-dimensional sphere) and we can probably finish the proof similarly as in \cite{HM}.


\vskip 10pt

In \cite{HM} we have in fact first recalled that each homeomorphism (defined on connected set) is either sense-preserving or sense-reversing and then we have showed that for $p>[n/2]$ each sense-preserving homeomorphism must satisfy $J_f\geq 0$ (and sense-reversing must satisfy $J_f\leq 0$). Naturally we wanted to claim that the topological notion of orientation (sense-preserving or sense-reversing) and analytical one ($J_f\geq 0$ or $J_f\leq 0$) are the same but unfortunately it turned out to be more difficult by the the following result of Hencl \cite{H}. 

\begin{thm}\label{jaczero}
Let $1\leq p<n$. There is a homeomorphism $f\in W^{1,p}([0,1]^n,[0,1]^n)$ with $f(x)=x$ for $x\in\partial [0,1]^n$ such that $J_f(x)=0$ almost everywhere. 
\end{thm}

It follows that such a mapping has many unintuitive properties, e.g. there is a set $N\subset (0,1)^n$ with
\eqn{qqq}
$$
|N|=0\text{ and }|f(N)|=1\text{ and thus also }|(0,1)^n\setminus N|=1\text{ and }|f((0,1)^n\setminus N)|=0. 
$$
For more results in this direction see e.g. \v{C}ern\'y \cite{C}, D'onofrio, Hencl and Schiattarella \cite{DHS}, Liu and Mal\'y \cite{LM}  and Oliva \cite{O}.

\open{Can the Jacobian of $W^{2,1}$ (or $WBV$) homeomorphism vanish} Let $n\geq 3$ and $1\leq q<n/2$. Is it possible to construct a Sobolev $W^{2,q}$ homeomorphism with 
$J_f=0$ a.e.? The same problem could be asked also for $f\in WBV$, i.e. for homeomorphism with first order derivative in the space of bounded variation. 

By Sobolev embedding theorem we know that for $q\geq n/2$ we have $W^{2,q}\subset W^{1,n}$ and hence $J_f\equiv 0$ is impossible as homeomorphisms in $W^{1,n}$ satisfy the Lusin $(N)$ condition (they map null sets to null sets) and thus \eqref{qqq} is impossible. The counterexample from \cite{H} for $n=2$ is constructed in a way that on almost every line parallel to coordinate axe we have derivative like  
$
\left( \begin{array}{cccc}
1&0\\
0&0 \\
\end{array} \right)
$
on a dense set and also derivative like
$
\left( \begin{array}{cccc}
0&0\\
0&1 \\
\end{array} \right)
$
on a dense sets. It follows that the variation of $Df$ on this line is clearly infinite and thus $Df$ fails the ACL condition and we definitely do not have $f\in W^{2,1}$ (or even $f\in WBV$). When I first thought about this problem it seemed to me that this cannot be difficult and there must be a simple proof but I failed to find any proof so far. 

\vskip 10pt


\end{document}